\documentclass[12pt]{amsart}   
\title[Euler-Lagrange Equations] {On derivation of Euler-Lagrange Equations for incompressible energy-minimizers}

\author[Nirmalendu Chaudhuri]{Nirmalendu Chaudhuri}
\address{ School of Mathematics and Applied Statistics University
of Wollongong, Wollongong, NSW 2522, Australia}
\email{chaudhur@uow.edu.au}

\author[Aram L. Karakhanyan]{Aram L. Karakhanyan}
\address{Department of Mathematics, University of Texas at Austin,Austin, Texas 78712, USA}
\email{aram@math.utexas.edu}

\thanks{$2000$ {\it Mathematics Subject Classification.\/} Primary
  35J60,  42A40,  73C50, 73V25.}

\date{July 24, 2008}
\usepackage{color}
\usepackage{graphics}
\usepackage[dvips]{graphicx}
\usepackage{wrapfig}
\usepackage{amsmath,amssymb}
\usepackage{amsfonts}
\usepackage[english]{babel}
\usepackage{enumerate}
\usepackage{graphicx}
\usepackage{mathrsfs}
\usepackage{mathrsfs}
\usepackage{color}
\usepackage[dvips]{graphicx}

\theoremstyle{plain}
\newtheorem{theorem}{Theorem}[section]
\newtheorem{lemma}[theorem]{Lemma}

\theoremstyle{remark}
\newtheorem{remark}[theorem]{Remark}
\theoremstyle{definition}

\renewcommand\epsilon\varepsilon 
\renewcommand\div{\operatorname{div}} 
\newcommand\supp{\operatorname{supp}} 

\numberwithin{equation}{section}
\newcommand\e{{\epsilon}}

\newcount\outlineon
\newcommand{\outline}[1]{\ifnum\outlineon>0\
\\\noindent\fbox{\begin{minipage}{10cm}{\footnotesize
#1}\end{minipage}}\\\\\fi}
\newcommand{\exclude}[1]{}

\newcommand{\h}[1]{\textbf{#1}}
\newcommand{\cal} {\mathcal}
\newcommand{\pa} {\partial}
\newcommand{\al} {\alpha}

\newcommand{\noi} {\noindent}

\newcommand{\ep} {\epsilon}
\newcommand{\na} {\nabla}
\newcommand{\be} {\begin{equation}}
\newcommand{\ee} {\end{equation}}
\newcommand{\la} {\lambda}
\newcommand{\de} {\delta}
\newcommand{\Om} {\Omega}

\def\lan{\langle}
\def\ran{\rangle}

\newcommand{\M} {\mathbb{M}}

\newcommand{\R} {\mathbb{R}}

\renewcommand\epsilon\varepsilon 
\renewcommand\div{\operatorname{div}} 
\numberwithin{equation}{section}
\numberwithin{equation}{section} \DeclareFontFamily{OT1}{pzc}{}
\DeclareFontShape{OT1}{pzc}{m}{it}%
             {<-> s * [1.250] pzcmi7t}{}
\DeclareMathAlphabet{\mathpzc}{OT1}{pzc}%
                                 {m}{it}

\def\XXint#1#2#3{{\setbox0=\hbox{$#1{#2#3}{\int}$}
     \vcenter{\hbox{$#2#3$}}\kern-.5\wd0}}

\makeindex

\begin{document}




\maketitle



\begin{abstract}
 We prove that any distribution $q$ satisfying the equation $\na q=\div\h f$ for some 
 tensor $\h f=(f^i_j),\,f^i_j\in h^r(U)$ ($1\leq r<\infty$) -the {\it local Hardy space}, $q$ is in $h^r$, and is locally represented by the sum of singular integrals of $f^i_j$ with Calder\'on-Zygmund kernel. As a consequence, we prove the existence and the local representation of the hydrostatic pressure $p$ (modulo constant) associated with incompressible elastic energy-minimizing deformation $\h u$ satisfying $|\na\h u|^2,\,|{\rm cof}\na\h u|^2\in h^1$. We also derive the system of Euler-Lagrange equations for incompressible local minimizers $\h u$ that are in  the space $K^{1,3}_{\rm loc}$ (defined in (\ref{k})); partially resolving a long standing problem. For H\"older continuous pressure $p$, we obtain partial regularity  of area-preserving minimizers. 
 \end{abstract}

\section{Introduction}

Let $\Om\subset\R^n$, $n\geq2$ be a bounded Lipschitz material
body. For Mooney-Rivlin or Neo-Hookean materials
\cite{Ball}, \cite{TO}, \cite{Og}, such as vulcanized rubber, in the
equilibrium state, one is interested in minimizing the elastic
energy
\begin{equation}\label{var1}
E[\h w]:=\int_{\Om}L(\na \h w(x))dx\,,
\end{equation}
for incompressible $W^{1,2}$-deformations $\h
w:\Om\subset\R^n\to\Bbb R^n$, subject to its own boundary condition,
and corresponding to a given smooth bulk energy $L:\M^{n\times n}\to
\Bbb R$.  Let us define the subspace $K^{1,r}$ for $1\leq
r<\infty$, by
\begin{equation}
K^{1,r}(\Om,\R^n):=\left\{\h w\in W^{1,r}(\Om,\R^n): {\rm cof}\,\na {\h w}\in
L^r(\Om,\mathbb M^{n\times n})\right\}\,,
\label{k}
\end{equation}
where $W^{1,r}$ denotes the usual {\it Sobolev spaces} (see for example, \cite[Chapter 7]{GT}) and ${\rm cof}\,P$ is the {\it cofactor} matrix,  whose $ij$-th entry is
$(-1)^{i+j}$ times the determinant of $(n-1)\times(n-1)$ submatrix
obtained by deleting the $i$-th row and the $j$-th column from the
$n\times n$ matrix $P$. Using the identity $P^t{\rm cof}P=Id_n\det P$, it follows that 
$\det\na\h w\in L^1$ for any $\h w\in K^{1,2}$. Since $|P|=|{\rm cof}\,P|$ for any
$P\in\M^{2\times 2}$, the function spaces $K^{1,r}$ and 
$W^{1,r}$ are equal in $\R^2$. Let us denote the admissible set of
deformations
\begin{equation}
\cal A:=\left\{\h w\in K^{1,2}(\Om,\R^n):   \det\na\h w =1\,\,\,{\rm a.e.}\,\,{\rm in}\,\,{\rm \Om}\right\},
\label{t1} \end{equation}
  We call ${\h u}\in\cal A$ to be a  {\it local minimizer} of
$E[\cdot]$ if and only if
\begin{equation}
E[{\h u}]\leq E[{\h w}]\quad {\rm for}\,\,{\rm all}\,\,{\h w}\in\cal
A\,\,\, {\rm and}\,\,\,{\rm supp}\,({\h w}-{\h u})\subset\Om\,.
\label{t2} \end{equation}
 Under the hypothesis that the energy density $L$ is smooth, {\it polyconvex} (convex function of minors) \cite{Ball}
 and satisfies the growth condition
 \begin{equation}
 C_1(|X|^2+|{\rm cof}\,X|^2)-C_2\leq L(X)\leq C_3(1+|X|^2+|{\rm cof}\,X|^2),
 \label{mrc}
 \end{equation}
 for all $X\in\M^{n\times n}$, for some $C_1>0$, $C_2\geq 0$, $C_3>0$, where $|X|^2:={\rm trace}(X^tX)$, using direct methods in the calculus of variations together with weak
 continuity of the determinant, J. Ball \cite{Ball} proved the existence of local
 minimizers $\h u\in\mathcal A$ of the energy $E[\cdot]$. An example
of polyconvex $L$ satisfying the growth condition (\ref{mrc}) is the stored-energy for incompressible isotropic Mooney-Rivlin materials in $\R^3$, given
by
\begin{equation}
L(X)=\frac{\mu_1}{2}(I_1(X)-3)+\frac{\mu_2}{2}(I_2(X)-3)\,,
\label{mr}
\end{equation}
where $I_1(X):={\rm trace}(C)=|X|^2$, $ I_2(X):=\frac{1}{2}
\left[\big({\rm trace}(C)\big)^2-{\rm trace}(C^2)\right]=|{\rm
cof}\,X|^2$, are the first two principle invariants of the right
Cauchy-Green strain tensor $C:=X^tX$ and $\mu_1$, $\mu_2$ are
positive material constants.

\medskip

Though the existence of the local minimizers of $E[\cdot]$ in $\mathcal
A$ is known for over 30 years, the existence of integrable hydrostatic pressure
associated with such minimizers, the derivation of system of
Euler-Lagrange equations, and the partial regularity for such
minimizers remains a challenging open problem. In this article we
prove the following results:
\medskip

\begin{enumerate}

\item The $h^r$ ($1\leq r<\infty$) -integrability and local representation of any distribution $q$ satisfying the equation $\na q=\h f$, where $\h f:=(f^i_j)$, $f^i_j\in h^r$, the {\it local $r$-Hardy spaces}. (Theorem
\ref{th1}) \bigskip

\item The existence of a pressure $p\in L^{r}_{\rm loc}$
if the minimizer is $\h u\in K^{1,2r}_{\rm loc}$ for some $r>1$.
(Theorem \ref{th2}) \bigskip

\item The existence of a pressure $p\in h^1$ if the
minimizer $\h u$ satisfies the conditions $|\na\h u|^2,\,|{\rm cof}\na\h u|^2\in h^1$. (Theorem
\ref{th2})
\bigskip

\item The validity of the Euler-Lagrange equations if the minimizer
is $\h u\in K^{1,3}_{\rm loc}$. (Theorem
\ref{th3}). The pair $(\h u,p)$ satisfies the system 
 \begin{equation} {\rm div} \,\left[DL(\na
\h u(x))- p(x)\,{\rm cof}\,(\na \h u(x))\right]=0\quad\quad{\rm
in}\;\;\Om, \label{el1}
\end{equation}
where the divergence is taken in each rows. 
\medskip

\item The partial regularity of $W^{1,3}$ area-preserving
minimizers $\h u$ for which the hydrostatic pressure $p$ is H\"older
continuous with exponent $0<\alpha<1$. (Theorem \ref{th4})
\medskip
\end{enumerate}

The $L^2$-version of the result in ({\bf I}) is classical (see,
\cite[Remark 1.4, p 11]{Te} ), and plays an important role in
incompressible fluids \cite{Te}. The result in ({\bf I}) is a
crucial ingredient in proving ({\bf II}) \& ({\bf III}). The
$h^1 $-version of ({\bf I}) is quite delicate and to the
best of our knowledge, it is new and may be of
independent interest. For the case $r>1$, it follows that $\na q\in W^{-1, r}$, and
adapting the classical functional-analytic approach demonstrated for
$r=2$ (see \cite{Te}, \cite{TO}), or arguing directly by duality, and solving
the equation of the type
$$
\div\h w=f\quad {\rm in}\,\,\, V\subset\subset U,\,\,\,\h w=0\quad {\rm in}\,\,\,\pa V,
$$
\cite[p. 472-474]{Ev2}, one can prove that $q\in L^{r}_{\rm loc}$.
However, both of these approaches fail to give informations for the critical case $r=1$ and 
does not give a representation of $q$. Whereas, our unified singular integral approach is self-contained, simple and provides the local $h^r$-estimate, as well as the local representation of $q$. The main ideas in our proof is to represent the localized-mollified distribution $q$ in terms of the Newtonian potential
in $\R^n$ and finding its uniform bound in $h^r$, by using Calder\'on-Zygmund estimate \cite{FS},  \cite{CZ}. Finally we show that the local representation of $q$ consists the sum of Calder\'on-Zygmund type singular integrals of the tensor $\h f$ (see equation (\ref{sq4}) in Section 4).
\medskip

For the case $n=2$, under the stronger hypothesis that the local
minimizers of $E[\cdot]$ are classical ($C^{1,\al}$-diffeomorphism), 
namely in the Sobolev space $W^{2,r}$ for some $r>2$, LeTallec and Oden \cite{TO} established
the system of equations in (\ref{el1}). For $n=2$, Bauman, Owen and Phillips \cite{BOP3} 
proved that if a minimizer is in $W^{2,r}$ for some $r>2$, then it is smooth. 
For such $W^{2,r}$, $r>2$ minimizers, the authors in \cite{BOP3} argued directly on the level
of the Euler-Lagrange equations exploring the existence of integrable hydrostatic pressure. Evans and Gariepy \cite{EG} proved that any {\it non-degenerate}, Lipschitz area-preserving local minimizers of $E[\cdot]$ are $C^{1,\al}(\Om_0)$, for some $0<\al<1$ for a dense open subset $\Om_0\subset \Om$. We believe that the Euler-Lagrange equations (\ref{el1}) 
that we derived for $K^{1,3}$-minimizers may be useful in understanding the partial regularity of such minimizers, as evidenced by the result in ({\bf V}). 
\medskip

In order to prove the existence of an integrable pressure $p$ associated with
the local minimizer $\h u$, we only require the additional mild
assumption $|\na \h u|^2\log(2+|\na u|^2)$, $|{\rm cof}\na\h u|^2\log(2+|{\rm cof}\na\h u|^2)\in L^1_{\rm loc}$.  For $n=2$, to derive the system of equilibrium equations (\ref{el1}) for $(\h u,p)$ in $\Om$,
we need $\h u$ to be in $W^{1,3}$, whereas the best-known previous result in this direction were for $W^{2,r}$-minimizers for some $r>2$. 
\medskip

We organize the paper as follows. In Section 2 we prove ({\bf I});
in Section 3 we prove ({\bf II}) \& ({\bf III}); in Section 4 we
prove ({\bf IV}), and finally in Section 5 we prove ({\bf V}).
Throughout this article $C$ is a generic absolute constant depending
on $n$, $U$, $\Om$, $\h u(\Om)$, $V\subset\subset\h u(\Om)$, $r$, and $L$. 
Its value can vary from line to line, but each line is valid with $C$ being a pure positive number.

\section{Local integrability of solutions $\na q=\div\h f$}

We recall some of the basic definitions and terminologies of Hardy spaces. Let $1\leq r<\infty$. A distribution $f$ belongs to $H^r(\R^n)$ if and only if $f\in L^r(\R^n)$ and $R_j(f)\in L^r(\R^n)$ (see for example, \cite[Proposition 3, p. 123]{St3}) for $j=1,\cdots,n$, where $R_j$ is the Riesz transform of $f$ given by 
$$
R_j(f)(x):=\lim_{\ep\to 0}c_n\int_{|y|\geq\e}\frac{y_j}{|y|^{n+1}} f(x-y)\,dy,\,\, c_n:=\frac{\Gamma\left(\frac{n+1}{2}\right)}{\pi^{(n+1)/2}},$$
so that $\widehat {R_j(f)}(\xi)=i\frac{\xi_j}{|\xi|}\widehat f$. In short, we will write $H^r(\R^n)$ as simply $H^r$. For $f\in H^r$, the norm is defined as 
$$\|f\|_{H^r}:=\|f\|_{L^r}+\sum_{j=1}^n\|R_j(f)\|_{L^r}.$$ 
A standard result \cite[p. 237]{St2} states that a positive function $f$, the Riesz transform $R_jf\in L^1_{\rm loc}$ if and only if $f\log(2+f)\in L^1_{\rm loc}$. For $1<r<\infty$, a classical result asserts that $f\in H^r$ if and only if $f\in L^r$, see \cite[p. 220]{St2}. The celebrated Fefferman duality theorem \cite{Fe}, \cite[Theorem 2]{FS}, \cite[Theorem 1, p. 142]{St3} asserts that the dual of $H^1$ is the BMO, the functions of bounded mean oscillations. The following theorem is due to Calder\'on-Zygmond \cite{CZ}, Stein \cite[Theorem 3, p. 39]{St2}, and Stein-Fefferman \cite[Corollary 1, p. 149-151]{FS}. 

\begin{theorem}[\bf Calder\'on-Zygmond, Fefferman-Stein] Let $1\leq r<\infty$ and $f\in H^r$. Let $G$ be a $C^1$ function on $\R^n\setminus\{0\}$ homogeneous of degree
$0$ with mean value $0$ over the unit sphere $\mathbb S^{n-1}$, that
is
\begin{equation} 
 \int_{\mathbb S^{n-1}}G(x)\,d\sigma(x)=0. \label{czg}
\end{equation}
Then the function defined as
\begin{equation}
T_0f(x):=\lim_{\de\to 0}\int_{|y|\geq
\de}\frac{G(y)}{|y|^n}\,f(x-y)\,dy \label{czf}
\end{equation}
exists a.e. and furthermore,
\begin{equation}
\|T_0f\|_{H^r}\leq C_{n,r}\|f\|_{H^r}\,.
 \end{equation} 
  \label{cza}
   \end{theorem}
   In particular, $R_j$'s are bounded linear operator on $H^r$, for any $1\leq r<\infty$. 
Let us recall the definition of {\it local Hardy spaces} introduced by Goldberg \cite{Go}.  A distribution $f$ on $\R^n$ is said to be in the local $r$-Hardy space, written as $f\in h^r$, if and only if the maximal function 
$$
\mathcal M_{\rm loc}f(x):=\sup_{0<\e<1}|(\rho_\e*f)(x)|
$$
is in $L^r$, where $\rho_\e:=\e^{-n}\rho(x/\e)$, is a standard approximation of the identity. The $h^r$ norm of $f$ is defined to be the $L^r$ norm of the maximal function $\mathcal M_{\rm loc}f$.  It follows that if $f\in h^r$ then $\eta f\in h^r$ for any smooth cut-off function and $H^r\subset h^r$.  For bounded Lipschitz domain $\Om\subset\R^n$, we adopt the definition of Hardy spaces $h^r(\Om)$ introduced by Miyachi \cite{Mi}. A distribution $f$ on $\Om$ is said to be in $h^r(\Om)$ if $f$ is the restriction to $\Om$ of a distribution $F$ in $h^r(\R^n)$, i.e., 
\begin{align}
\nonumber
h^r(\Om)&:=\{f\in\mathcal D'(\Om):\,\exists\, F\in h^r(\R^n),\,\,{\rm such\,\,that}\,\,F\big|_{\Om}=f\}\\\nonumber
&=h^r(\R^n)\big/\{F\in h^r(\R^n)\,:\,F=0\,\,{\rm on}\,\,\Om\}.
\end{align}
The norm on this space is the quotient norm: the infimum of $h^r$ norms of all possible extensions of $f$ in $\R^n$. For $1<r<\infty$ the spaces $h^r(\Om)$ is equivalent to $L^r(\Om)$. For smooth bounded domains $\Om$, the Theorem \ref{cza} is valid for $f\in h^1(\Om)$, see \cite{Mi}, \cite{CKS}. 
\medskip


\begin{theorem} Let $U\subset\R^n$, $n\geq2$ be a bounded Lipschitz domain and $1\leq r<\infty$. Let $\h f=(f^i_j)$ such that $f^i_j\in h^r(U)$, for $1\leq i,j\leq n$. Then the distribution $q:C^{\infty}_0(U)\to\R$ defined by
\begin{equation}
\na q={\rm div} \,\h f\iff \lan \na q,\h v\ran=-\int_{U}\h f(x):\na
\h v(x)\,dx\, \label{diveq}
\end{equation}
for $ \h v\in C^{\infty}_0(U,\R^n)$, is in $h^r(V)$, for any $V\subset\subset U$ where $A:B:={\rm
trace}(A^tB)=\sum_{ij}a^i_{j}b^i_{j}$, for $A,B\in\M^{n\times n}$.
Furthermore, $q$ is locally represented by sum of singular integrals
of $f^i_j$ (see equation (\ref{sq4})), and  for any
$V\subset\subset U$, there exists $C>0$, depending only on $U$, $V$
and $r$ such that
$$
\|q\|_{h^r(V)}\leq C\|\h f\|_{h^r(V)}\,.
$$
\label{th1}
\end{theorem}


\noi{\bf Proof of Theorem \ref{th1}}. Let $U\subset\R^n$, $n\geq2$
be a Lipschitz domain. Let $\h f:=(f^i_j)\in \M^{n\times n}$ and
$f^i_j\in h^r(U)$, for $1\leq r<\infty$ and
$1\leq i,j\leq n$. Let $q\in\mathcal D'(U)$, such that
\begin{equation}
\na q={\rm div} \,\h f\quad{\rm in}\,\,\,\, \mathcal D'(U)\,.
\label{diveq1}
\end{equation}
 Our idea is to mollify the equations in (\ref{diveq1}) and obtain
uniform bound for the mollified $q$, by using Calder\'on-Zygmund
estimate. Let $V\subset\subset U$ be a sub-domain and $0<\e<{\rm
dist}(V,\pa U)$. Let $\rho_\e$ be the usual mollification kernel,
and define convolution  $q_{\e}:V\to\R$ by
$$
q_\e(x)=(q*\rho_\e)(x):=\lan q,(\rho_\e)_x\ran\quad {\rm for}\,\,\,
x\in V, \quad{\rm where}\,\,\, (\rho_{\ep})_x(y)
:=\rho_{\ep}(y-x),\,\,\,y\in U
$$
Then by the standard properties of the mollification
\cite[Proposition 1, p492]{DL}, $q_{\ep}$ is smooth and for any
$1\leq i\leq n$
$$
\frac{\pa}{\pa {x_i}}(q*\rho_\e)=\frac{\pa q}{\pa {x_i}}*\rho_{\e}=q* \frac{\pa\rho_\e}{\pa {x_i}}.
$$
Hence mollifying the system of equations in (\ref{diveq1}), we obtain
\begin{equation}
\na q_{\ep}={\rm div}\,\h f_{\ep}\quad {\rm in}\,\,\, V,
\label{gq}
\end{equation}
where the divergence is taken in each rows of $\h
f_{\ep}:=\bigg((f^i_j)_\e\bigg)$, and $(f^i_j)_\e:=f^i_j*\rho_{\ep}$
is the mollification of $\h f$.  Since $f^i_j\in h^r(U)$, we conclude that
\begin{equation}
\label{stro}
(f^i_j)_{\e}\to f^i_j\quad {\rm strongly}\,\, {\rm in}\,\,\,
h^r(V)\quad {\rm as}\,\,\ep\to 0,
\end{equation}
 for all $1\leq i,j\leq n$. Applying the divergence operator to
 the both sides of the above equation, we
obtain
\begin{eqnarray}
\Delta q_\e=\div(\div\h f_\e) \quad{\rm in}\,\,\,V.
\label{lap}
\end{eqnarray}
Since there is no control on the boundary values, we need to
localize the equation (\ref{lap}). Let $W\subset\subset
V\subset\subset U$. Let $\eta\in C^{\infty}_0(\R^n)$, $0\leq \eta\leq 1$ be a cut-off
function such that $ \eta\equiv 1\ \textrm{in}\, W $ and
$\eta\equiv0$ outside $V$. Let $\bar q_\e:=\eta q_\e$ be the
localization of $q_{\ep}$. Then $\bar q_\e$ is the solution of
Poisson equation
\begin{equation}
   \Delta \bar q_\e=\bar f_\e\quad {\rm in}\,\,\, \R^n,
\end{equation}
where
\begin{align}
\label{fe}
\bar f_\e &:=\eta\Delta q_\e+2\lan\na q_\e,
\na\eta\ran+q_\e\Delta\eta\\\nonumber
&=\eta\div(\div\h f_\e)+2\lan \div\h f_{\ep},\na\eta\ran +q_\e\Delta\eta.
\end{align}
Therefore $\bar q_\e$ is represented by the Newtonian
potential of in $\R^n$. In other words,
\begin{eqnarray}\label{loc-pres}
  \bar{q}_\e(x)=-\int_{\R^n}\Phi(x-y)\bar f_\e(y)\,dy\,,
\end{eqnarray}
where $\Phi$ is fundamental solution of the Laplace equation in $\R^n$ and is given by
\begin{equation}
\Phi(x):=\left\{\begin{array}{ll}
-\frac{1}{2\pi}\log |x| & {\rm if}\,\,\,\, n=2\\\\
\frac{1}{n(n-2)\alpha(n)}\,\frac{1}{|x|^{n-2}} & {\rm if}\,\,\,\, n\geq 3\,,
\end{array}\right.
\label{green}
\end{equation}
for $x\in\R^n\setminus\{0\}$, and
$\alpha(n):=\frac{\pi^{n/2}}{\Gamma(\frac{n}{2}+1)}$ is the volume
of the unit ball in $\R^n$. Using (\ref{fe}) in (\ref{loc-pres}), we
obtain
\begin{eqnarray}
\bar q_\e(x)&=&-\int_{\R^n}\eta(y)\Phi(x-y)\div(\div\h f_\e)\,dy\\\nonumber
&&+2\int_{\R^n}\big(\lan \div\h f_{\ep},\na\eta\ran
+q_\e\Delta\eta\big)\Phi(x-y)\,dy\\\nonumber
   &:=&- I^1_\e(x)-2I^2_{\ep}(x)-I^3_\e(x)\,,
   \end{eqnarray}
where
\begin{eqnarray}
\label{ae1}
  I^1_{\ep}(x)&:=& \int_{\R^n}\eta(y)\,\Phi(x-y)\,\div(\div\h f_\e(y))  \,  dy\\
  \label{ae2}
  I_{\ep}^2(x)&:=& \int_{\R^n}\left\lan\div \h f_\e(y), \na\eta(y)\right\ran\,\Phi(x-y)\,dy\\
  \label{ae3}
  I_{\ep}^3(x)&:=&\int_{\R^n}q_\e(y)\,\Delta\eta(y)\,\Phi(x-y)\,dy
\end{eqnarray}
By direct computations, observe that, for $1\leq i,j\leq n$
\begin{eqnarray}
\label{be1} \big(\eta\,\Phi\big)_{y_i}&=&\eta_{y_i}\Phi(y)
-\frac{1}{\omega_n}\,\frac{\eta\,y_i}{|y|^n}\,,\\
\label{be2} \big(\eta \,\Phi\big)_{y_iy_j}
&=&\eta_{y_iy_j}\,\Phi(y)-\frac{1}{\omega_n}\frac{y_i\eta_{y_j}+y_j\eta_{y_i}}{|y|^n}\\\nonumber
&&-\frac{1}{\omega_n}\,\left(\delta_{ij}-n\,\frac{y_iy_j}{|y|^2}\right)\frac{\eta}{|y|^n},
\end{eqnarray}
where $\delta_{ij}$ is the Kr\"onecker delta and
$\omega_n:=n\alpha_n$ is the surface area of the unit sphere
$\mathbb S^{n-1}$. We now establish an uniform local $h^r$-estimates ($1\leq r<\infty$)
for $q_{\ep}$ through the following steps.
\smallskip

\noi
{\bf Step 1: Limit of $I^3_{\ep}$}.\, Let us fix $x\in  W\subset\subset V\subset\subset U$.
Since $\Delta\eta=0$ on $W$, the integrand in $I^3_\e(x)$ is smooth.
Since $q_\e$ is determined up to a constant, we can
add a constant to $y\mapsto\Delta\eta(y)\Phi|x-y|$, if nessecary, to
ensure that it has vanishing integral. For each fixed $x\in W$, let $\h v_x:V\to\R^n$ be the
solution of the Dirichlet problem
\begin{eqnarray}
  \left\{\begin{array}{ll}
    \div \h v_x(y)=\Delta\eta(y)\Phi(x-y)  \quad & \mbox{for}\  y\in V\\
    \h v_x= 0 & \mbox{on }\  \partial V\,.
  \end{array}
  \right.
  \label{psd}
\end{eqnarray}
Then using (\ref{psd}), integrating by parts, and the convergence of
$\h f_\e$ in (\ref{ae3}), we obtain
\begin{eqnarray} \label{i30}
I^3_{\ep}(x)&=&  \int_{\R^n}q_\e(y)\Delta\eta(y)\Phi(x-y)\,dy\\\nonumber
&=&\int_{\R^n}q_\e(y)\div\h v_x(y)\,dy\\\nonumber
  &=&-\int_{\R^n}\lan\na q_\e(y),\h v_x(y)\ran\,dx\\\nonumber
&=& -\int_{\R^n}\lan\div\h f_\e(y),\h v_x(y)\ran\,dy\\\nonumber &=&
\int_{\R^n}\h f_\e(y): \na_y \h v_x(y)\ran\,dy\\\nonumber
 &\rightarrow&\int_{\R^n}\h f(y):\na_y \h v_x(y)\,dy\quad{\rm as}\,\,\e\to 0\\\nonumber
&:=&I^3_0(x)\,\quad{\rm for}\,\,\,x\in W\subset\subset V.
\end{eqnarray}
Since $\h f_\e\to\h f$ strongly in $h^r(V,\M^{n\times n})$, it follows that $I^3_\e\to I^3_0$
strongly in $h^r(W)$.
\medskip

\noindent {\bf Step 2: Limit of $I^2_{\ep}$}. Let us fix $x\in
W\subset\subset V\subset\subset U$. Integrating by parts, invoking
(\ref{be1}) and letting $\e\to 0$ we have
\begin{eqnarray}
\label{i20} I^2_{\ep}(x) &=&\int_{\R^n}\bigg\lan\div\h
f_{\ep}(y),\Phi(x-y) \,\na\eta(y)\bigg\ran\,dy\\\nonumber
&=&-\int_{\R^n}\h
{f}_{\ep}:\na_y\bigg(\Phi(x-y)\,\na\eta\bigg)\,dy\\\nonumber
&=& -\int_{\R^n}\h f_{\ep}:\bigg(\Phi(x-y)\,\na^2\eta
  -\frac{(y-x)\otimes\na\eta}{\omega_n\,|y-x|^n}\bigg)\,dy
\\\nonumber
&\rightarrow&\,
 -\int_{\R^n}\h f:\bigg(\Phi(x-y)\,\na^2\eta
  -\frac{(y-x)\otimes\na\eta}{\omega_n\,|y-x|^n}\bigg)\,dy\\\nonumber
 &:=&I^2_0(x)\quad x\in W\,.
 \end{eqnarray}
  Using the strong convergence of $\h f_\e$ in $h^r(V)$, again it follows that $I^2_\e\to I^2_0$ in $h^r(W)$.
 \medskip

{\bf Step 3: Limit of $I^1_{\ep}$}. Integrating by parts twice
the integral in (\ref{ae1}) and using (\ref{be2})
\setlength\arraycolsep{1pt}
\begin{eqnarray}\nonumber
I^1_{\ep}(x)&=&\int_{\R^n}\div \div\h f_\e(y)\,\,\eta(y)\Phi(x-y)\,dy
\\\nonumber
&=&\int_{\R^n}\h f_\e(y):\na^2_y\bigg(\eta(y)\Phi(x-y)\bigg)\,dy
\\\nonumber
&=&\int_{\R^n}\h f_\e(y):\left(\Phi(x-y)\,\na^2\eta(y)
-\frac{1}{\omega_n}\frac{\na\eta\otimes(y-x)+(y-x)\otimes\na\eta}{|x-y|^n}\right)\,dy\\\nonumber
&&-\frac{1}{\omega_n}\int_{\R^n}\h f_\e(y):
\left(Id_n-n\frac{(y-x)\otimes(y-x)}{|x-y|^2}\right)\frac{\eta}{|x-y|^n}\bigg)\,dy\\\nonumber
  &:=&I^{11}_\e(x)+I^{12}_\e(x),\quad x\in W,
\end{eqnarray}
where $Id_n$ is the $n\times n$ identity matrix. Using the
convergence of $\h f_\e$, observe that as $\e \to 0$,
\begin{eqnarray}\nonumber \label{xyt}
 I^{11}_\e(x)&:=&\int_{\R^n}\h f_\e:\left(\Phi(x-y)\,\na^2\eta
-\frac{\na\eta\otimes(y-x)+(y-x)\otimes\na\eta}{\omega_n\,|x-y|^n}\right)dy
\\\nonumber &\rightarrow&
\int_{\R^n}\h f:\left(\Phi(x-y)\,\na^2\eta
-\frac{\na\eta\otimes(y-x)+(y-x)\otimes\na\eta}{\omega_n\,|x-y|^n}\right)dy
\\&:=&I^{11}_0(x),\,\quad x\in W.
\end{eqnarray}
In order to estimate $I^{12}_\e$, define the kernels
$\Omega_{ij}:\R^n\setminus\{0\}\to\R$ by
\begin{equation}
\Omega_{ij}(y):=\de_{ij}-n\,\frac{y_iy_j}{|y|^2},\quad
y\in\R^n\setminus\{0\},\,\,\, \,\,\,i,j=1,\cdots,n.
\end{equation}
Since $n\alpha_n=\omega_n$, integrating by parts, observe that for
any $i,j=1,\cdots,n$,
\begin{align}\nonumber
\int_{\mathbb S^{n-1}}\Omega_{ij}(y)\,d\sigma(y)&=\int_{\mathbb
{S}^{n-1}}(\de_{ij}-ny_iy_j)\,d\sigma(y)\\\nonumber
&=\omega_n\de_{ij}-n\int_{\mathbb S^{n-1}}y_iy_j\,d\sigma(y)\\\nonumber
&=\omega_n\de_{ij}-n\int_{B_1}\frac{\pa}{\pa y_j}y_i\,dy\\\nonumber
&=\omega_n\de_{ij}-n\de_{ij}\alpha_n\\\nonumber &=0\,.
\end{align}
Hence each $\Omega_{ij}$ satisfies all the conditions of
Calder\'on-Zygmund Kernel \cite{St2}. Therefore,
\begin{equation}
I^{12}_{\ep}(x):=-\frac{1}{\omega_n}\int_{\R^n}\eta\h f_\e:
\left(Id_n-n\frac{(y-x)\otimes(y-x)}{|x-y|^2}\right)\frac{dy}{|x-y|^n}\bigg)
 \label{s1}
\end{equation}
is the sum of Calder\'on-Zygmund type singular integrals
with the homogeneous kernel $\Omega_{ij}$. Since $\h f\in h^r(U,\M^{n\times n})$, $1\leq r<\infty$, by Theorem \ref{cza} $I^{12}\in h^r(W)$. Furthermore, the following
sum of singular integrals
\begin{equation}
\label{s2}
I^{12}_0(x):=-\frac{1}{\omega_n}\int_{\R^n}\eta\h f:
\left(Id_n-n\frac{(y-x)\otimes(y-x)}{|x-y|^2}\right)\frac{dy}{|x-y|^n}
\end{equation}
exists for almost every $x\in W\subset\subset V$ and is in $h^r(W)$.
From the singular integrals (\ref{s1}) and (\ref{s2}), by Theorem \ref{cza}, we have
$$
I^{12}_{\e}(x)-I^{12}_0(x)=-\frac{1}{\omega_n}\sum_{i,j=1}^n\int_{\R^n}\bigg(\eta(f^i_{j})_{\ep}(y)-
\eta f^i_{j}(y)\bigg) \frac{\Omega_{ij}(x-y)}{|x-y|^n}\,dy\,. \label{calz}
$$
Hence there exists $C:=C(V,W,r)>0$ such that 
\begin{equation}
\| I^{12}_{\ep}-I^{12}_0\|_{h^r(W)}\leq C\sum_{j=1}^n \|(f^i_j)_{\ep}-f^{i}_j\|_{h^r(V)}\to 0\quad{\rm as}\,\,\,\e\to 0.
\label{fa}
\end{equation}

\noi{\bf Step 4: Explicit representation of $q$}.\, To complete the
proof, let us define the potential $q:W\to\R$ by
$$
q(x):=-\big(I^{11}_0(x)+I^{12}_0(x)+2I^2_0(x)+I^3_0(x)\big)\,.
$$
Then from (\ref{i30}), (\ref{i20}), (\ref{xyt}), and (\ref{fa}), we conclude that $q_{\ep}\rightarrow q$ strongly in $h^r_{\rm loc}$ for any $1\leq r<\infty$, and hence and $q$ is
represented as
\begin{eqnarray}
\label{sq4} q(x)& =&\int_{U}\h f:\left(\Phi(x-y)\na^2\eta-\na_y\h
v_x\right)dy\\\nonumber &&+ \frac{1}{\omega_n}\int_{U}\h
f:\bigg(\na\eta\otimes(y-x)-(y-x)\otimes\na\eta\bigg)\frac{dy}{|x-y|^n}\\\nonumber
&&+\frac{1}{\omega_n}\int_{U}\eta\h f:
\left(Id_n-n\frac{(y-x)\otimes(y-x)}{|x-y|^2}\right)\frac{dy}{|x-y|^n}\end{eqnarray}
for any $x\in W.$ Since $q$ is the strong limit of the family $q_{\ep}$ in $W$, it is
independent of the choice of the cut-off function $\eta$. This completes the proof of Theorem 1.1. \qed

\section{First Variation of Energy and the existence of hydrostatic pressure}
  Let $\Om\subset\R^n$, $n\geq 2$ be a smooth, simply connected and bounded domain and let $L:\M^{n\times n}\to\R$ be smooth function. We are now in a position to establish the existence of integrable hydrostatic pressure associated with 
 the local minimizers of the energy 
\begin{equation}\label{var9}
E[\h w]:=\int_{\Om}L(\na \h w(x))dx\,,
\end{equation}
for incompressible $W^{1,2}$-deformations $\h
w:\Om\subset\R^n\to\Bbb R^n$. By direct computation, observe that Mooney-Rivlin bulk-energy 
given by
\begin{equation}
L(X)=\frac{\mu_1}{2}(|\na\h u|^2-3)+\frac{\mu_2}{2}(|{\rm cof}\na\h u|^2-3)\,,
\label{mr9}
\end{equation}
satisfies the following. 

\begin{eqnarray}\nonumber
DL&=&\mu_1 P+\mu_2 
\left(\begin{array}{ccc}
{\rm cof}(SQ)^1_{1}:(SP)^1_{1} & -{\rm cof}(SQ)^1_{2}:(SQ)^1_{2} 
& {\rm cof}(SQ)^1_{3}:(SP)^1_{3}\\\\
-{\rm cof}(SQ)^2_{1}:(SP)^2_{1} & {\rm cof}(SQ)^2_{2}:(SP)^2_{2}
& -{\rm cof}(SQ)^2_{3}:(SP)^2_{3}\\\\
{\rm cof}(SQ)^3_{1}:(SP)^3_{1} & -{\rm cof}(SQ)^3_{2}:(SP)^3_{2}
& {\rm cof}(SQ)^3_{3}:(SP)^3_{3}
\end{array} \right),
\end{eqnarray}
where $Q:={\rm cof}P$, and $(SX)^i_{j}$ is the $2\times 2$ submatrix obtained by deleting the $i$-th row and the $j$-th column of the matrix $X\in M^{3\times 3}$. Furthermore, the Cauchy-Green strain tensor is given by 
\begin{eqnarray}\nonumber
(DL(P))^tP&=&\mu_1P^tP+\mu_2
 \left(\begin{array}{ccc}
\big|Q_2\big|^2+\big|Q_3\big|^2& -\big\lan Q_1,Q_2\big\ran 
&-\big\lan Q_1,Q_3\big\ran \\\\
-\big\lan Q_1,Q_2\big\ran &\big|Q_1\big|^2+\big|Q_3\big|^2
&-\big\lan Q_2,Q_3\big\ran \\\\
-\big\lan Q_1,Q_2\big\ran &-\big\lan Q_2,Q_3\big\ran
&\big|Q_1\big|^2+\big|Q_2\big|^2\end{array} \right)
\end{eqnarray}
for all $P\in \M^{3\times3}$, where $Q_i:=({\rm cof}P)_i:=\big(({\rm cof}P)^i_1,({\rm cof}P)^i_2,({\rm cof}P)^i_3\big)$ be the $i$-th row of ${\rm cof}P$, $i=1,2,3$.  Motivated by the above calculations, assume that $L$ satisfies the 
following growth condition.
\begin{equation}
\max\bigg( |L(P)|,\big|(DL(P))^tP\big|\bigg)\leq C\big(1+|P|^2+|{\rm cof}P|^2\big),\label{qe}
\end{equation}
 for some $C>0$, for any $P\in \M^{n\times n}$. 
 \medskip
 
 Now we prove the existence of an integrable hydrostatic pressure $q$ on the deformed domain $\h u(\Om)$ and establish an explicit representation of the pressure $q$ in terms of Calder\'on-Zygmund type singular integrals of the Cauchy-Green strain $\tilde\sigma:=(DL(\na\h u))^t\na\h u)\circ \h u^{-1}$ in $\h u(\Om)$. Our proof consists of deriving the first variation of the energy $E[\cdot]$, obtaining the equation $\na q=\div\tilde\sigma$, and then finally use Theorem \ref{th1}.

\begin{theorem} Let $L:\M^{n\times n}\to\R$ be smooth and satisfies the growth condition (\ref{qe}). 
Assume that ${\h u}\in \mathcal A$ be a continuous and injective
local minimizer of $E[\cdot]$, such that $|\na\h {u}|^2, |{\rm cof}\na\h u|^2\in
h^r_{\rm loc}(\Om)$ for some $1\leq r<\infty$. Then there
exists a scalar function $q\in h^r_{\rm loc}(\h{u}(\Om))$, such that
$$
 \|q\|_{h^{r}(V)}\leq C\left(\big\||\na\h u|^2\big\|_{h^r(\h u^{-1}(V))}+\big\||{\rm cof}\na\h u|^2\big\|_{h^r(\h u^{-1}(V))}\right) ,\quad
 V\subset\subset \h u(\Om),$$
 for some $C>0$ (depending on $r$, $V$, $n$ and $\h u(\Om)$) and the pair $(\h u,q)$ satisfies the integral identity
\begin{equation}
\int_{\Om}DL(\na \h u(x))\,:\,\na (\h v\circ \h u)\,dx= \int_{\h u
(\Om)}q(y)\,\div \h v(y) \,dy
  \label{imd}
\end{equation}
for all $\h v \in C_0^\infty(\h u(\Om), \R^n)$, where $A:B:=\rm
{tr}(A^tB)=\sum_{i,j=1}^n a^i_{j}b^j_{j}$ is the scalar product on
$\M^{n\times n}$.
 \label{th2}
\end{theorem}

\begin{remark}  Let $W\subset\subset V\subset\subset\h u(\Om)$, and $\eta\in
C^{\infty}_0(V)$ be a cut-off function such that $\eta\equiv 1$ on
$W$. Then $q$ is locally represented as
\begin{eqnarray}
\label{sqi} q(x)& =&\int_{V}
\tilde\sigma:\left(\Phi(x-y)\na^2\eta-\na_y\h
v_x\right)dy\\\nonumber &&+\frac{1}{\omega_n} \int_{V}\tilde\sigma
:\bigg(\na\eta\otimes(y-x)-(y-x)\otimes\na\eta\bigg)\frac{dy}{|x-y|^n}\\\nonumber
&&+\frac{1}{\omega_n}\int_{V}\eta\tilde\sigma:
\left(Id_n-n\frac{(y-x)\otimes(y-x)}{|x-y|^2}\right)\frac{dy}{|x-y|^n}
,\nonumber
\end{eqnarray}
for any $x\in W$, where $\Phi$ is Newtonian potential in $\R^n$
defined in (\ref{green}) and $\h v_x$ as defined in (\ref{psd}). 
\end{remark}

\begin{remark}
In the study of regularity of finite energy deformations, \v{S}ver\'{a}k  \cite{Sv 1}
proved that for any $W^{1,n}$-deformation $\h w$ with ${\rm det}\, \na\h w(x)>0$,
a.e., there exists a continuous function $\omega$ on $\R$ with
$\omega(0)=0$ such that
$$
|\h w(x)-\h w(y)|\leq \omega(|x-y|),\quad{\rm for}\,\,\,{\rm any}\,\,\,x,y\in\Omega\subset\subset\R^n.
$$
It is also well-known any $W^{1,n}$-deformation $\h w$ for which 
 the {\it distortion} function $K(\cdot,\h w):=|\na\h w(\cdot)|^n/{\rm det}\,\na\h
{w}(\cdot)\in L^r$ for some $r>n-1$, is a homeomorphism. Thus in
particular, area-preserving $W^{1,r}$ ($r>2$)-deformations in the plane are
continuous and open maps. However, in general for $n\geq 3$, any deformation 
$\h w\in K^{1,2}$ may be totally discontinuous, see \cite[p. 119]{Sv 1}. 
\end{remark}

 In order to prove Theorem \ref{th2}, we establish the 
following first variation of the energy integral $E[\cdot]$.

\begin{lemma} {\bf First Variation.} 
 Let $\h u\in{\mathcal A}$ be a local minimizer of $E[\cdot]$. We further assume
that $\h u$ is a continuous and an injective map. Then
$\h u$ satisfies the following integral identity
\begin{equation}
 \int_{\Om}  DL(\na \h  u(x)):\na(\h v\circ \h u)(x)\,dx=0\,,
\label{MainId}
\end{equation}
for all smooth, compactly supported and divergence free vector fields $\h
v$ on $\h u(\Om)$. \label{wd}
\end{lemma}

\noindent{\bf Proof:} By the invariance of domain $\h u(\Om)$ is open and $\h u:\Om\to\h u(\Om)$ is a homeomorphism. Let ${\bf v}\in C_0^\infty(\h u(\Om), \R^n)$ be a vector field with $\div{\bf v }=0$. For each $y\in \h u(\Om)$,
 consider the unique smooth flow $\phi(y,\cdot):\R\to \h u(\Om)$
  given by
  \begin{equation}
    \frac{d\phi}{dt}(y,t)={\bf v}(\phi(y,t)) \quad \textrm{in}\quad \R,\quad\phi(y,0)=y.
  \end{equation}
  Using the relations ${\displaystyle \frac{\pa}{\pa P^i_j}{\rm det}\,P=({\rm cof}\,P)^i_j}$ and
  $P\,({\rm cof}\,P)^t=Id_n\,{\det}\,P$, by a direct calculations we
  observe that

\setlength\arraycolsep{2pt}
\begin{equation}
  \frac{d}{dt}\left(\det\na_y
  \phi(y,t)\right)=\det\na_y\phi(y,t)\div{\bf
  {v}}=0.
  \label{z}
  \end{equation}
Since $\det\na_y \phi(y,0)=1$, from (\ref{z}) it
follows that $\det\na_y\phi(y,t)=1$ for all $t\in\R$ and $y\in \h u(\Om)$.
Consider the map $\h w:\Om\times\R\to \h u(\Om)$ defined by
$$\h w(x,t):=\phi(\cdot,t)\circ {\bf u}\,(x)=\phi(\h u(x),t)\quad
{\rm for}\,\,\,{\rm any}\,\,\, t\in\R,\,\,\,x\in\Om.
$$
Let $V:=\supp {\bf v}\subset \h u(\Om)$, then ${\bf v}(\h u(x))=0$
for ${\bf u}(x)\not\in V$. This in conjunction with the uniqueness
of $\phi$ implies that $\phi({\bf u}(x),t)={\bf u}(x)$ for all
points $x$ such that ${\bf u}(x)\not \in V$. Since $\Om$ is bounded,
$\bf u$ is continuous and $V$ is compact, $\Om'={\bf u}^{-1}(V)$ is
a compact subset of $\Om$. Hence $\supp({\bf w}(x,t)-{\bf
u}(x))\subset \Om'$. Furthermore, $\det\na_x{\h w(x,t)}=\det\na_y
\phi(y,t)\,\det\na {\h u(x)}=1$. Therefore, ${\bf w}(\cdot,t)\in
\cal A$ and $\supp (\h u-\h w(\cdot,t))\subset\Om$ for all $t\in\R$.
Since $\bf u$ is a local minimizer of $E[\cdot]$,
$$
E[{\bf u}]\leq E[{\h w(\cdot,t)}]\quad\quad{\rm for}\;{\rm all}\quad
t\in\R.
$$
Thus in particular,
\begin{align}
0&=\left.\frac{d}{dt}\int_{\Om}L(\na \h w(x,t))\,dx\right|_{t=0}\nonumber\\
&=\sum_{i,j=1}^{2}\left.\int_{\Om}L^i_j(\na \h w(x,t))\,
\frac{d}{dt}\left(\frac{\pa w^i}{\pa x_j}(x,t)\right)dx\right|_{t=0}\nonumber\\
&=\sum_{i,j=1}^{2}\left.\int_{\Om}L^i_j(\na \h w(x,t))\,
\frac{\pa}{\pa x_j}\left(\frac{d\phi^i}{dt}(\h u(x),t)\right)
dx\right|_{t=0}\nonumber\\
&=\sum_{i,j=1}^{2}\left.\int_{\Om}L^i_j(\na \h w(x,t))\,
\frac{\pa}{\pa x_j}\left(v^i(\phi(u(x),t)\right)dx\right|_{t=0}\nonumber\\
&=\sum_{i,j=1}^{2}\int_{\Om}L^i_j(\na \h u(x))\,
\frac{\pa}{\pa x_j}\left(v^i(\h u(x))\right)dx\nonumber\\
&=\int_{\Om} DL(\na \h u(x))\,:\,\na(\h v\circ \h
u)(x)\,dx\,,\nonumber
\end{align}
for all smooth, compactly supported and divergence free vector
fields on $\h u(\Om)$, where ${\displaystyle L^i_j(P):=\frac{\pa
L}{\pa p^i_{j}}(P)}$ . This proves the Theorem. \qed
\bigskip

{\bf Proof of Theorem \ref{th2}:} 
Let $1\leq r<\infty$ and $U'\subset\subset U$. Let $\h u\in\mathcal A$ be a local minimizer of $E[\cdot]$ such that $|\na\h u|^2\in h^r$ and $|{\rm cof}\,\na\h u|^2 \in h^r(U')$ for some $1\leq r<\infty$. Assume further that $\h u:\Om\to \h u(\Om)$ is continuous and bijective map. 
\medskip

Now define $\h g=(g^1,\cdots,g^n)\,:C^1_0(\h u(\Om),\R^n)\to\R$ by
\begin{equation}
\lan \h g,\h v\ran:=
\int_{\Om}  DL(\na \h  u(x)):\na(\h v\circ \h u)(x)\,dx,
\label{dis}
\end{equation}
for all $\h v=(v^1,\cdots,v^n) \in C^1_0(\h u(\Om),\R^n)$. In view of the volume constraint and growth condition (\ref{qe}), it follows that
\begin{equation}
|\lan \h g,{\h v}\ran|  \leq C\big(1+\|\na u\|_{L^2(\Om)}+\|{\rm cof}\na u\|_{L^2(\Om)}
\big)\,\|\na{\h
v}\|_{L^{\infty}({\h u}(\Om))}, \label{des}
\end{equation}
for any ${\h v}\in C^{1}_0(\h {u}(\Om),\R^n)$. Hence $\h g$ is a
continuous linear functional on $C^1_0({\h u}(\Om),\R^n)$. Using the the first variation 
(\ref{MainId}), we conclude that 
\begin{equation}
\lan \h g,{\h v}\ran=0\quad\forall\,\,{\h v}\in C^{1}_0(\h {u}(\Om),\R^n),\,\,{\rm div}\,\h v=0\,.
\label{div0}
\end{equation}
Hence there exists $ q \in \mathcal D'(\h u(\Om))$ ( see \cite[Proposition 1.1, p10]{Te}),
such that \begin{equation}
\h g =-\na q\,\quad{\rm in}\,\,\mathcal D'(\h u(\Om),\R^n)\label{diq}
\end{equation}
modulo translation of a constant. In order to obtain $h^r$ estimates of $q$, for $1\leq i,j\leq n$, let us define $\sigma^i_{j}:\Om\to\R$ by
\begin{eqnarray}\label{sigma}
\sigma^i_{j}(x):=\sum_{k=1}^n
  L^i_{k}(\nabla \h u(x))\,\frac{\partial u^j}{\partial x_k}(x)\quad{\rm for}\,\,
  x\in\Om,
  \end{eqnarray}
  so that, the Cauchy-Green strain tensor on $\Om$  is given by 
  \begin{equation}
  \sigma:=\big(\sigma^i_j\big)=\big(DL(\na\h u))^t\na\h u
  \label{sip}
  \end{equation}
   Define the $ij$-th component of the Cauchy-Green Strain tensor $\tilde\sigma^i_{j}$ on the deformed domain $\h u(\Om)$ by 
\begin{equation}
\tilde\sigma^i_{j}:=\sigma^i_{j}\circ \h u^{-1}\quad{\rm on}\,\,\h u(\Om),\,\,\, i,j=1,\cdots,n.
\label{tildesig}
\end{equation}
The growth condition $|\sigma^i_j|\leq C(|\na\h u|^2+|{\rm cof}\na\h u|^2)$ and $|\na\h u|^2,|{\rm cof}\na\h u|^2\in L\log L$ yields $\tilde \sigma^i_j\in h^1(V)$. If $\h u\in K^{1,2r}_{\rm loc} (\Om,\R^n)$, $1<r<\infty$, from the definition of $\sigma^i_{j}$, $\tilde\sigma^i_{j}$,  and the
condition (\ref{qe}) on $L$, it follows that
\begin{eqnarray}\label{l1}
\int_{V}|(\tilde\sigma^i_{j}|^{r} 
 &=&\int_{\h u^{-1}(V)}|\sigma^i_{j}|^{r}\\\nonumber
& \leq &C\left( \|\nabla
  \h u\|^{2r}_{L^{2r}(\h u^{-1}(V))}+  \|{\rm cof} \nabla
  \h u\|^{2r}_{L^{2r}(\h u^{-1}(V))} \right)
  \,,\end{eqnarray}
for any $V\subset\subset \h u(\Om)$. Therefore, if $|\na\h u|^2\in h^r$ and $|{\rm cof}\,\na\h u|^2 \in h^r_{\rm loc}$ for some $1\leq r<\infty$, from (\ref{l1}), we have 
$$\sigma:=\big(\sigma^i_{j}\big)\in
h^{r}_{{\rm loc}}(\Om,\M^{n\times n})\quad{\rm and}\,\,\, 
\tilde\sigma:=\big(\tilde\sigma^i_{j}\big)\in h^{r}_{\rm
loc}(\h u(\Om),\M^{n\times n}).
$$
 Observe that, from the
definition of $\h g$ in (\ref{dis}), $\sigma^i_{j}$ in
(\ref{sigma}), $\tilde\sigma^i_{j}$ in (\ref{tildesig}), and change
of variables,
\begin{align}\label{G-e}
 \big\lan \h g,\h v\big\ran&=\sum_{i,k=1}^n  \int_{\Om}
  L^i_{k}(\nabla \h u(x))\,\frac{\partial }{\partial x_k}(v^i\circ \h u)(x)\,dx\\\nonumber
 &= \sum_{i,j,k=1}^n  \int_{\Om}
  L^i_{k}(\nabla \h u(x))\,\frac{\partial v^i }{\partial y_j}(\h u(x))\frac{\pa u^j}{\pa x_k}(x)\,dx\\\nonumber
 &= \sum_{i,j=1}^n  \int_{\Om}\sigma^i_j(x)\,\frac{\partial v^i }{\partial y_j}(\h u(x))\,dx\\\nonumber
&= \int_{\Om}  \sigma(x):\na_{\h u}\h v(\h u(x))
dx\\\nonumber
&=\int_{\h u(\Om)}  \tilde\sigma(y):\na\h v(y)\,
dy\\\nonumber
   &=-\big\lan\div\tilde\sigma,\h v\big\ran
  \end{align}
for any $v\in C^1_0(\h u(\Om),\R^n)$. Hence
\begin{equation}
\h g =-\div\tilde\sigma\,\quad{\rm in}\,\,\mathcal D'(\h u(\Om),\M^{n\times n})
\label{gtil}
\end{equation}
where the divergence is taken in each rows. Therefore, combining (\ref{diq}) and (\ref{gtil}), we get
\begin{equation}
\na q =\div\tilde\sigma\,\quad{\rm in}\,\,\mathcal D'(\h u(\Om),\M^{n\times n}).
\label{qsig}
\end{equation}
By taking $\h f=\tilde \sigma$, and $U=V\subset\subset\h u(\Om)$ in (\ref{qsig}), from Theorem \ref{th1}, we conclude that $q\in h^{r}_{\rm loc}(\h u(\Om))$, it satisfies the local representation (\ref{sqi}), and
 \begin{eqnarray}
 \|q\|_{h^{r}(V)}&\leq &C\|\tilde\sigma\|_{h^{r}(V)}\\\nonumber
 &\leq& C\left(\||\na\h u|^2\|_{h^r(\h u^{-1}(V))}+
 \||{\rm cof}\na\h u|^2\|_{h^r(\h u^{-1}(V))}\right) ,
 \label{qest}
 \end{eqnarray}
 for any $V\subset\subset \h u(\Om)$, for some $C>0$, depending on $r$, $V$, $n$ and $\h u(\Om)$. Since $q\in L^1_{\rm loc}$, from (\ref{diq}), it follows that
$$
\lan \h g,\h v\ran=-\lan \na q,\h v\ran=\lan q,\div\h v\ran=\int_{\h u(\Om)}q(y)\,\div \h v(y)\,dy.
$$
for any $\h v\in C^1_0(\h u(\Om),\R^n)$. Hence
\begin{equation}
\int_{\Om} DL(\na \h u(x)):\nabla(\h v\circ \h u)(x)dx
=\int_{\h u(\Om)} q(y)\,{\rm div}\,\h v(y)\,dy,
\label{rep}
\end{equation}
for any $ \h v\in C^1_0(\h u(\Om),\R^n)$. This proves the Theorem. \qed

\section{Derivation of Euler-Lagrange Equations}

\begin{theorem} Let $\Om\subset\R^n$, $n\geq 2$, be a smooth,
simply connected and bounded domain. Let ${\h u}\in \mathcal
A\cap K^{1,s}_{\rm loc}(\Om,\R^n)$ for some $s\geq 3$ be a continuous and injective 
local minimizer of $E[\cdot]$. Then the hydrostatic pressure $p:=q\circ \h
u\in L^{s/2}_{\rm loc}(\Om)$, and  the pair $(\h u,p)$ satisfies
\begin{equation}
\int_{\Om}DL(\na \h u(x))\,:\,\na \phi(x)\,dx= \int_{\Om}p(x)\,{\rm
cof}\,(\na \h u(x))\,:\,\na {\phi(x)}\,dx\,, \label{eli}
\end{equation}
for all $\phi \in C^{1}_0(\Om,\R^n)$, where $q \in L^{s/2}_{\rm
loc}(\h u(\Om))$ as in Theorem \ref{th2}. In other words, the pair $(\h
u,p)$ satisfies the system of Euler-Lagrange equations
$$
{\rm div} \,\left[DL(\na \h u(x))- p(x)\,{\rm cof}\,(\na \h
u(x))\right]=0\quad\quad{\rm in}\;\;\Om,
$$
in the sense of distribution, where the divergence is taken in each
rows. \label{th3}
\end{theorem}
\medskip

 {\bf Proof. } Let $\Om\subset\R^n$ be a smooth, simply connected domain.
 Recall that $K^{1,s}:=\{\h w\in W^{1,s}\,:\,{\rm cof}\,\na\h w\in
L^s\}$ and $\mathcal A:=\{{\bf w}\in K^{1,2}(\Om,\R^n): \,
\det\na\bf w =1\ \textrm{a.e.}\}$. Let $\h u\in \mathcal A\cap
K^{1,s}_{\rm loc}(\Om,\R^n)$, $s\geq 3$ be a continuous injective
local minimizer of the functional $E[\cdot]$. By Theorem \ref{th2}, there
exists $q\in L^{s/2}_{\rm loc}$ such that the pair $(\h u, q)$
satisfies the identity (\ref{rep}).  Let $\h u^{-1}:\h
u(\Om)\to\Om$ be the inverse of $\h u$. Then using the
volume-constraint we obtain
$$
\na_y\h u^{-1}(y)=(\na_x\h u(x))^{-1}=\big({\rm cof}\,\na\h u(x)\big)^t,\quad y=\h u(x),
$$
and hence by the change of variables
$$
\int_{\h u(\Om)}|\na\h u^{-1}(y)|^2dy=\int_{\Om}|{\rm cof}\,\na\h u(x)|^2dx<\infty.
$$
 Using the relation ${\rm cof}\, (XY)={\rm cof}\, X\,{\rm cof}\, Y$, for $X,Y\in\M^{n\times n}$, observe that
$$
Id_n={\rm cof}\,\big(\na_y\h u^{-1}\, \na\h u\big)={\rm cof}\big(\na_y\h u^{-1}\big)\,{\rm cof}\big(\na\h u\big)
={\rm cof}\big(\na_{y}\h u^{-1}\big)\,(\na \h u)^{-t},
$$
and hence
$${\rm cof}\big(\na\h u^{-1}\big)=(\na\h u)^t\,.
$$
Since $\h u\in K^{1,s}_{\rm loc}(\Om,\R^n)$, it follows that ${\h u}^{-1}\in
K^{1,s}_{\rm loc}(\h u(\Om),\Om)$ for $s\geq 3$. Let
$V\subset\subset \h u(\Om)$ and $\phi \in C^{1}_0(\h
u^{-1}(V),\R^n)$. Then the composition
  $\phi\circ {\h u}^{-1}\in W^{1,s}_0(V,\R^n)$.
  Hence there exists ${\h v}_{\e}\in C^1_0(V,\R^n)$
 such that ${\h v}_{\e}\to \psi:=\phi\circ {\h u}^{-1}$ strongly in
$W^{1,s}(V,\R^n)$ as $\e\to 0$. Let $U:=\h u^{-1}(V)$. Then H\"older inequality yields
\begin{eqnarray}\nonumber
  \int_{U} DL(\na \h u):\bigg(\na(\h v_\e\circ \h u)
  &-&\na(\psi\circ \h u)\bigg)dx\\\nonumber
  &=&\int_{U}(\na \h u)^t DL(\na \h u):
  \bigg(\na_z\h v_\e(\h u)-\na_z\psi(\h u)\bigg)dx\\\nonumber
  &\leq& C\|\na \h u\|_{L^{2s'}(U)}\,\|\na(\h v_\e-\psi) \|_{L^s(V)},
\end{eqnarray}
where $s':=s/(s-1)$. Notice that $s\geq 3$ yields $2s'\leq s$ and
hence $\na\h u\in L^{s}_{\rm loc}(\Om)\subseteq  L^{2s'}_{\rm
loc}(\Om)$. Therefore, from (\ref{dis}) we obtain
\begin{align}\label{dcm}
\lan \h g,\h v_{\e}\ran&=\int_{\h u^{-1}(V)} DL(\na \h
u(x)):\nabla(\h v_{\e}\circ \h u)(x)\,dx\\\nonumber &\to\int_{\h
u^{-1}(V)} DL(\na \h u(x)):\nabla(\phi\,\circ\,\h u^{-1}\circ \h
u)(x)\,dx \quad{\rm as}\,\,\,\e\to 0\\\nonumber &=\int_{\h
u^{-1}(V)} DL(\na \h u(x)):\nabla\phi(x)\,dx\,.
\end{align}
Since $ \na\h u,\,{\rm cof}\,\na\h u\in L^s_{\rm loc}$, $q\in L^{s/2}_{\rm loc}$ and
$L^{s/2}_{\rm{loc}}\subseteq L^{s/(s-1)}_{\rm loc}$ for $s\geq 3$,
applying change of variables in (\ref{rep}), and letting $\e\to 0$ we obtain
   \begin{align}\label{sem}
\left \lan \h g,{\h v_{\e}}\right\ran
   &= \int_{V}q(y)\,\,{\rm trace}\left({\na \h v_{\e}}(y)\right)\,dy\,.\\\nonumber
  &= \int_{\h u^{-1}(V)}q(\h u(x))\,\,{\rm trace}\bigg({\na_{\h u} \h v_{\e}}(\h u(x))\bigg)\,dy
  \\\nonumber
&=\int_{\h u^{-1}(V)}q({\h u}(x))\,\,{\rm trace}\bigg (\na({\h
v_{\e}}\circ{\h u})(x)\,\,({\rm cof}\,\na\h
u(x))^t\bigg)\,dx\\\nonumber
  &= \int_{\h u^{-1}(V)}q({\h u}(x))\,\,{\rm cof}\,(\na {\h
  u}(x)):\na({\h v_{\e}}\circ{\h u})(x)\,dx,\\\nonumber
  &\to\int_{\h u^{-1}(V)}q({\h u}(x))\,\,{\rm cof}\,(\na {\h
  u}(x)):\na(\phi\,\circ\,\h u^{-1}\circ{\h u})(x)\,dx\\\nonumber
   &=\int_{\h u^{-1}(V)}q({\h u}(x))\,\,{\rm cof}\,(\na {\h
  u}(x)):\na\phi(x)\,dx\,.
    \end{align}
  Hence from (\ref{dcm}) and (\ref{sem}) we obtain
 $$
\int_{\h u^{-1}(V)} DL(\na \h u(x)):\nabla\phi(x)\,dx =\int_{\h
u^{-1}(V)}q({\h u}(x))\,\,{\rm cof}\,(\na {\h
  u}(x)):\na\phi(x)\,dx\,,
$$
   for any $\phi\in C^1_0(\h u^{-1}(V),\R^n)$. Finally choose a sequence of
 smooth, simply connected sets $V_k\subset\subset
V_{k+1}\subset\subset \h u(\Om)$ sub-domains such that $\h
{u}(\Om)=\cup_{k=1}^{\infty}V_k$. Utilizing the foregoing arguments, there exists $q_k\in
L^{s/2}(V_k)$, $k\geq 1$ such that
\begin{equation}
\int_{\h u^{-1}(V_k)} DL(\na \h u):\nabla\phi=\int_{\h
u^{-1}(V_k)}q_k(\h u)\,\,{\rm cof}\,(\na {\h
  u}):\na\phi\,,
\label{intm}
\end{equation}
for $\phi\in C^1_0(\h{u}^{-1}(V_k),\R^n)$. Since $\h u$ is locally
volume-preserving homeomorphism, $\Om=\cup_{k=1}^{\infty}\h
u^{-1}(V_k)$ is an open covering of $\Om$ and
$\h{u}^{-1}(V_k)\subset\subset \h{u}^{-1}(V_{k+1})$. Using the identity
$\div\,{\rm cof}\,\na\h{u}(x)=0$ and invertibility of $\na \h u(x)$,
from (\ref{intm}) it follows that $q_k$ is unique up to a
translation of a constant. Thus adding constant terms as necessary
to each $q_k$, we deduce from (\ref{intm}) that for each fixed
$k\geq 1$
$$
q_i(z)=q_k(z)\quad{\rm for}\,\,z\in V_i,\quad 1\leq i\leq k.
$$
We finally define $q:\h u(\Om)\to\R$ as $q(z):=q_k(z)$, for $z\in
V_k$, so that $q\in L^{s/2}_{\rm loc}(\h u(\Om))$.  This proves that
for any $\phi\in C^1_0(\Om,\R^n)$, the pair $(\h u,q)$ satisfies
$$
\int_{\Om} DL(\na \h u(x)):\nabla\phi(x)\,dx=\int_{\Om}q(\h
u(x))\,\,{\rm cof}\,(\na {\h
  u}(x)):\na\phi(x)\,dx\,.
$$
Now let us define the pressure $p$ on $\Om$ by
$$p(x):=q(\h u(x))\quad{\rm for} \,\,x\in\Om.
$$
Then for any $k\geq 1$,
 $$
 \int_{\h u^{-1}(V_k)} |p(x)|^{s/2}=\int_{\h u^{-1}(V_k)}|q(\h u(x))|^{s/2}dx
 =\int_{V_k}|q(z)|^{s/2}dz<\infty,
 $$
 and hence $p\in L^{s/2}_{\rm loc}(\Om)$ and the pair $(\h u,p)$ satisfies
\begin{equation}
 \int_{\Om} DL(\na \h u(x)):\nabla\phi(x)\,dx=\int_{\Om}p(x)\,\,{\rm cof}\,(\na {\h
  u}(x)):\na\phi(x)\,dx\,,  \label{intf}
\end{equation}
for any $\phi\in C^1_0(\Om,\R^n).$ In other words, $(\h u,p)$
satisfies the system of Euler-Lagrange equations
$$
{\rm div}\left[DL(\na{\h u(x)})-p(x)\,{\rm cof}\,(\na {\h
u(x)})\right]=0,\quad {\rm in} \,\,\,\,\Om.
$$
in the sense of (\ref{intf}).  This completes the proof. \qed

\section{Partial Regularity of area-preserving minimizers}

For $n=2$, as a consequence of the Euler-Lagrange equations (\ref{el1}), together with the
standard elliptic estimates \cite{GM}, we establish the following theorem.
\medskip

 \begin{theorem}
   Let $\Om\subset\R^2$ be a smooth, bounded simply connected domain and  let $L:\M^{2\times 2}\to\R$ be smooth, uniformly convex, such that $DL$ has linear growth and $D^2L$ is bounded. Let $\h u\in W^{1,3}(\Om,\R^2)$ be an area-preserving minimizer of the energy $E[\cdot]$. Furthermore, assume that the associated hydrostatic pressure $q$ on the deformed domain $\h u(\Om)$ is $C^{0,\alpha}$ for some $0<\alpha<1$. Then $\na \h u$ are H\"older continuous on a dense open set $\Om_0\subset\Om$.    \label{th4}
 \end{theorem}

{\bf Proof.} Since $\h u\in W^{1,3}(\Om,\R^2)$ and $\h u$ is area-preserving, $\h u(\Om)$ is open and $\h u$ is a homeomorphism from $\Om$ to $\h u(\Om)$. By Theorem \ref{th3}, there exists $q\in L^{3/2}_{\rm loc}(\h u(\Om))$ and the pair 
$(\h u,q\circ \h u)$ satisfies the system 
\begin{equation}
\sum_{j=1}^2\frac{\pa}{\pa x_j}\left(\frac{\pa L}{\pa p^i_j}(\na\h u)-p(x)\left({\rm cof}\,\na {\h
u}\right)^i_j\right)=0,\quad {\rm in} \,\,\,\,\Om,\,\,\,i=1,2,
\label{2el}
\end{equation}
where $p:=q\circ \h u$.  Assume that $q\in C^{0,\al}(\h u(\Om))$. Since $\h u\in W^{1,3}$, Sobolev imbedding theorem yields 
$\h u\in C^{1/3}$, and hence $p(x)=q(\h u(x))$ is H\"older continuous with the exponent $\al/3$. 
Let $F:\Om\times \M^{2\times2}\to\R$ be the free-energy defined as 
$$
F(x,P):=L(P)-p(x)\det P\quad x\in\Om,\,\,\,P\in\M^{2\times 2},
$$
so that we can rewrite the nonlinear system (\ref{2el}) as 
\begin{equation}
\sum_{j=1}^2\frac{\pa}{\pa x_j}\left(A^i_j(x,\na\h u)\right)=0,\quad {\rm in} \,\,\,\,\Om,\,\,\,i=1,2,
\label{3el}
\end{equation}
where 
$$A^i_j(x,P):=\frac{\pa F}{\pa p^i_j}(x,P)=\frac{\pa L}{\pa p^i_j}(P)-p(x)({\rm cof}P)^i_j.$$ Let $U\subset\subset\Om$. Since $|{\rm cof}P|=|P|$ for any $P\in\M^{2\times 2}$, $|DL(P)|\leq C(1+|P|)$ and $D^2L(P)$ is bounded, 
\begin{equation}
|A^i_j(x,P)|\leq C(1+|P|),\,\,\,\left|\frac{\pa A^i_j}{\pa p^k_l}(x,P)\right|\leq C,
\label{con1}
\end{equation}
for any $x\in U,\,\,P\in\M^{2\times 2}$. By H\"older continuity of $p$, it follows that 
\begin{eqnarray}
\frac{|A^i_j(x,P)-A^i_j(y,P)|}{1+|P|}&=&|p(x)-p(y)|\frac{\left|({\rm cof}P)^i_j\right|}{1+|P|}\\\nonumber
&\leq& C |x-y|^{\al/3},
\end{eqnarray}
for any $x\in U,\,\,P\in\M^{2\times 2}$. By direct calculations and the ellipticity of $L$ it follows that 
\begin{eqnarray}
\frac{\pa A^i_j}{\pa p^k_l}(x,P)\xi_{ij}\xi_{kl}&=&\frac{\pa^2 F}{\pa p^i_jp^k_l}(x,P)\xi_{ij}\xi_{kl}\\\nonumber
&=& \frac{\pa^2 L}{\pa p^i_jp^k_l}(P)\xi_{ij}\xi_{kl}-2p(x)\det\xi\\\nonumber
&\geq& \la_0|\xi|^2-2p(x)\det\xi\\\nonumber
&:=&I(x,\xi),\quad{\rm for}\,\,\,P=(p^i_j),\,\xi=(\xi_{ij})\in\M^{2\times 2},
\label{ell}
\end{eqnarray}
where $\la_0>0$ is the ellipticity constant of $L$. Completing squares, observe that 
\begin{eqnarray}
  \label{it1}
    \frac{I(x,\xi)}{\la_0}&=& |\xi|^2-2\frac{p(x)}{\la_0}\det\xi\\\nonumber  
  & =&\xi^2_{11}+\xi^2_{12}+\xi^2_{21}+\xi^2_{22}-2\frac{p}{\la_0}\left(\xi_{11}\xi_{22}-\xi_{12}\xi_{21}\right)\\\nonumber
  &=&\left(\xi_{11}-\frac{p}{\la_0}\xi_{22}\right)^2+\left(\xi_{12}-\frac{p}{\la_0}\xi_{21}\right)^2\\\nonumber
  &&+\left(1-\frac{p^2}{\la_0^2}\right)(\xi_{22}^2+\xi_{21}^2).\\\nonumber
\end{eqnarray}
Similarly, we obtain 
\begin{equation}
  \frac{I(x,\xi)}{\la_0}= \left(\xi_{22}-\frac{p}{\la_0}\xi_{11}\right)^2+\left(\xi_{21}-\frac{p}{\la_0}\xi_{12}\right)^2+\left(1-\frac{p^2}{\la_0^2}\right)(\xi_{11}^2+\xi_{12}^2)
  \label{it2}
\end{equation}
Adding the identities (\ref{it1}) and (\ref{it2}), we obtain 
\begin{eqnarray}
\label{it3}
 2\frac{ I}{\la_0}&=&
 \left(\xi_{11}-\frac{p}{\la_0}\xi_{22}\right)^2+\left(\xi_{12}-\frac{p}{\la_0}\xi_{21}\right)^2\\\nonumber &&+\left(\xi_{22}-\frac{p}{\la_0}\xi_{11}\right)^2
 +\left(\xi_{21}-\frac{p}{\la_0}\xi_{12}\right)^2 +\left(1-\frac{p^2}{\la_0^2}\right) 
 |\xi|^2\\\nonumber
  &\geq &\left(1-\frac{p^2}{\la_0^2}\right)|\xi|^2.
\end{eqnarray}
Thus from (\ref{ell}) and (\ref{it3}), it follows that the map $P\mapsto A(\cdot,P)$ is {\it strongly elliptic} if there exists $\mu_0>0$ such that 
$$\frac{\pa L^i_j}{\pa p^k_l}(x, P)\xi_{ij}\xi_{kl}\geq \frac{\la_0}{2}
\left(1-\frac{p^2}{\la_0^2}\right)|\xi|^2\geq
\mu_0|\xi|^2, \quad{\rm for}\,\,x\in\Om,\,P,\xi\in\M^{2\times2},$$
which is equivalent to assume that
\begin{equation}\label{small-p}
p^2\leq \la_0^2-2\lambda_0\mu_0\,\,\Longrightarrow (p-\mu_0)^2\leq (\la_0-\mu_0)^2.
\end{equation}
Since $p$ is defined up to addition of arbitrary constant, thus the inequality 
(\ref{small-p}) is satisfied in subdomain $U\subset\subset\Omega$ if and only if 
\begin{equation}\textrm{osc}_U \,p<\la_0.
\label{osc1}
\end{equation}
Since $p$ is H\"older continuous, the estimate (\ref{osc1}) holds for any subdomain $U\subset\Om$ with 
sufficiently small diameter. Hence $A(x,P)$ is strongly elliptic in $P$ for each $x\in U\subset\subset\Om$, for sufficiently small diameter. This proves that $A^i_j(x,P)$ satisfies all the conditions of Giaquinta-Modica in \cite{GM} on $U\subset\subset\Om$, with diameter of $U$ being small. Hence by \cite[Theorem 1]{GM}, we conclude that $\na \h u$ is H\"older continuous on a dense open subset $U_0$ of $U$. By standard covering arguments we conclude the proof.\qed 
\bigskip

\noi{\bf Acknowledgement} This work was initiated while both the
authors were at the Australian National University, which was
supported by Australian Research Council. The second author was
partially supported by the National Science Foundation.

\enddocument